Excellent Student Thesis

# Estimate of the Convergence Rate of Finite Element Solutions to Elliptic Equations of Second Order with Discontinuous Coefficients[1] [2]

Jin-chao Xu

**Abstract**

In this paper, we consider elliptic boundary value problems with discontinuous coefficients and obtain the asymptotic optimal error estimate $\|u-u_k\|_{1,\Omega} \leq Ch|\ln h|^{1/2} \|u\|_{2,\Omega_1+\Omega_2}$ for triangle linear elements.

## 1. Introduction

Finite element convergence theory has been well established for second order elliptic problems with appropriately smooth coefficients (cf. [1]). However, we frequently encounter problems with discontinuous coefficients in practical engineering applications. Accordingly, it is necessary to study the the finite element convergence for these problems. There has been some work on this topic in the one-dimensional case (see for example [5]), but fewer results are known for higher dimensional cases, which are of more theoretical and practical importance. A simple case when the jump interface is a polygonal line has been discussed by Professor Feng Kang in [3], though no results have been mentioned for more general cases. This paper is devoted to a discussion on the convergence of piecewise linear finite element approximations on triangular meshes for the Dirichlet problem with discontinuous coefficients in 2D. The study indicates that, the discontinuity in the coefficients has only mild influence to the convergence of the finite element approximation.

For simplicity, we consider the following model problem

$$\begin{cases} -\partial_x(B\partial_x u) - \partial_y(B\partial_y u) + \sigma u = f(x,y), \\ \left(B\dfrac{\partial u}{\partial n}\right)_{1,S} = \left(B\dfrac{\partial u}{\partial n}\right)_{2,S}, \\ u|_\Gamma = 0. \end{cases} \quad (1.1)$$

Here $\Omega$ is a bounded connected domain in 2D, and the boundary $\Gamma = \partial\Omega$ is piecewise smooth and convex. A piecewise smooth curve $S$ divides $\Omega$ into two subdomains $\Omega_1$

---

and $\Omega_2$. The coefficient $B = B(x, y) \in C^1(\overline{\Omega}_i)$, when restricted on $\Omega_i$, for $i = 1, 2$, and $\sigma \in L^\infty(\Omega)$ and $f \in L^2(\Omega)$. Moreover, we assume

$$B(x, y) \geq B_0 > 0 \quad \text{and} \quad \sigma \geq 0 \text{ on } \Omega. \tag{1.2}$$

We introduce the bilinear form:

$$a(u, v) = \iint_\Omega [B(\partial_x u \partial_x v + \partial_y u \partial_y v) + \sigma uv] dxdy. \tag{1.3}$$

It is easy to show the following variational principle:

> Solving equation (1.1) is equivalent to the variational problem: find $u \in H_0^1(\Omega)$ such that for any $v \in H_0^1(\Omega)$ it holds that

$$a(u, v) = (f, v). \tag{1.4}$$

The bilinear form $a(u, v)$ defined by (1.3) is continuous and coercive on $H_0^1(\Omega)$. Assume that the problem (1.4) admits a solution $u \in H_0^2(\Omega, S) = \{u \mid u \in H_0^1(\Omega), u \in H^2(\Omega_1), u \in H^2(\Omega_2)\}$. The finite element approximation, $u_h$ in a finite-dimensional subspace $S_h \subset H_0^1(\Omega)$, satisfies the basic error estimate:

$$\|u - u_h\|_{1,\Omega} \leq C \inf_{v \in S_h} \|u - v\|_{1,\Omega}. \tag{1.5}$$

Here and throughout the paper, $C$ denotes a generic constant independent of $h$, $u$, and $v$. In particular, let $u_I \in S_h$ be the interpolant of $u$, then we have

$$\|u - u_h\|_{1,\Omega} \leq C\|u - u_I\|_{1,\Omega}. \tag{1.6}$$

Now we consider a triangulation $\Omega_h \subset \Omega$ with a boundary $\Gamma_h$ whose vertices all lie on $\Gamma$. We assume that every triangle intersecting at $S$ has two vertices on $S$, every non-smooth point on $S$ is set to be a vertex, and every triangle contains a disk whose radius is $ch$, where $h$ is the maximal diameter for triangles in $\Omega_h$ (implying $\Omega_h$ is a quasi uniform triangulation). Let $S_h$ be the space of continuous and piecewise linear polynomials defined on $\Omega$ that vanish on $\Gamma_h$. Then, $S_h$ is a subspace, $S_h \subset H_0^1(\Omega)$. The elements that do not intersect the jump interface, $S$, are referred to as *regular* elements, and the other elements are *irregular* elements. The main result of this paper is: the finite element solution $u_h \in S_h$ and exact solution $u \in H^2(\Omega, S)$ of the variational problem (1.4) satisfy the following error estimate:

$$\|u - u_h\|_{1,\Omega_h} \leq Ch|\ln h|^{1/2}\|u\|_{2,\Omega_1+\Omega_2}. \tag{1.7}$$

## 2. Proof of the result

First, we give two lemmas, which are estimates similar to the Sobolev embedding theorem (cf. [6]), though sharper.

**Lemma 1.** *Let $\Omega$ be the aforementioned planar domain and $f \in L^2(\Omega)$. Denote $p = (x, y)$, $Q = (\zeta, \eta)$, and let $B(p, Q)$ be a bounded function of $p, Q$, which is continuous whenever $p \neq Q$. Define*

$$v(p) = \iint_\Omega \frac{B(p, Q)}{|p - Q|} f(Q) dQ.$$

*Then there exists a constant C, such that for any measurable set $D \subset \Omega$ it holds that*

$$\iint_D v^2(p)\, dp \leqslant \frac{C}{\varepsilon} |D|^{1-\varepsilon} \|f\|^2_{L^2(\Omega)},$$

*where $|D|$ is the measure of D, and $\varepsilon < 1$ is an arbitrary positive constant.*

*Proof.* Let $q = \frac{2}{\varepsilon} > 2$, then by the Hölder inequality (cf. [6]), we have

$$\iint_D |v|^q\, dp \leqslant C^q q^{\frac{q}{2}} \|f\|^{q-2}_{L^2(\Omega)} \iint_D dp \iint_\Omega f^2(Q) r^{-\frac{1}{2}}\, dQ.$$

Change the order of integration to obtain

$$\iint_D |v|^q\, dp \leqslant C^q q^{\frac{q}{2}} \|f\|^q_{L^2(\Omega)}.$$

Making use of the Hölder inequality again, we get

$$\iint_D v^2\, dp \leqslant |D|^{1-\frac{2}{q}} \|v\|_{L^q(D)} \leqslant \frac{C}{\varepsilon} |D|^{1-\varepsilon} \|f\|^2_{L^2(\Omega)}.$$

This completes the proof. □

By this lemma and the Sobolev integral identity (cf. [1]), we can obtain the lemma below.

**Lemma 2.** *Let $D, \Omega, \varepsilon, C$ follow the definitions in Lemma 1. It holds for any $v \in H^1(\Omega)$ that*

$$\iint_D v^2 dp \leqslant \frac{C}{\varepsilon} |D|^{1-\varepsilon} \|v\|^2_{1,\Omega}.$$

We now turn to the proof of (1.7).

**1. Estimate on regular elements.**

Naturally, the error estimate on regular elements can be obtained by the well-known Bramble-Hilbert lemma (cf. [1]). Here, we adopt a method based on the Taylor expansion with respect to varying the base point of expansion (see for example [4]), which will play a fundamental role in the following analysis. Let $K$ be any regular element, with the vertices $p_i = (x_i, y_i)$; let $p = (x, y)$ be a varying point, and $M_i = [x_i + (x - x_i)t, y_i + (y - y_i)t]$ for $i = 1, 2, 3$ and $0 \leqslant t \leqslant 1$. Without loss of generality, we assume $v \in C^2(\Omega)$ in the analysis below. Then, by the Taylor expansion with the integral form remainder, we have

$$u(p) - u(p_i) = (x - x_i)\partial_x u(p) + (y - y_i)\partial_y u(p) - \int_0^1 t\partial_t^2 u(M_i)dt.$$

Using the barycentric coordinates $L_i$, we can write the error of the linear interpolation function $u_I$ as:

$$u(p) - u_I(p) = -\sum_{i=1}^{3} L_i \int_0^1 t\partial_t^2 u(M_i)dt. \tag{2.1}$$

By the Hölder inequality, we get

$$\|u - u_I\|_{0,K}^2 \leqslant 3 \sum_{i=1}^{3} \iint_K \left|\int_0^1 t\partial_t^2 u(M_i)dt\right|^2 dxdy$$

$$\leqslant 3 \sum_{i=1}^{3} \iint_K \int_0^1 t^2 |\partial_t^2 u(M_i)|^2 dtdxdy.$$

Exchange the order of integration, and for any fixed $t, i$, introduce a change of variables $\zeta = x_i + (x - x_i)t$, $\eta = y_i + (y - y_i)t$, to transform the domain $K$ to a subdomain $K_{i,t} \subset K$. Then, we obtain

$$\|u - u_I\|_{0,K}^2 \leqslant Ch^4 |u|_{2,K}^2. \tag{2.2}$$

Moreover, from (2.1), we have

$$\partial_x(u - u_I) = -\sum_i \partial_x L_i \int_0^1 t\partial_t^2 u(M_i)dt - \sum_i L_i \int_0^1 t\partial_x \partial_t^2 u(M_i)dt.$$

Using integration by parts, the last term in the above equation becomes

$$\sum_i L_i \int_0^1 t\partial_x \partial_t^2 u(M_i)dt = \sum_i L_i \int_0^1 t[\partial_t^2 u_x(M_i)t + 2\partial_t u_x(M_i)]dt$$

$$= \sum L_i t^2 \partial_t u_x(M_i)|_0^1 = 0.$$

Note that $|\partial_x L_i| \leqslant C/h$, and by the same technique above, we obtain

$$\|\partial_x(u - u_I)\|_{0,K}^2 \leqslant Ch^2 |u|_{2,K}^2, \tag{2.3}$$

$$\|\partial_y(u - u_I)\|_{0,K}^2 \leqslant Ch^2 |u|_{2,K}^2. \tag{2.4}$$

Combing all three inequalities above, we have

$$\|u - u_I\|_{1,K}^2 \leqslant Ch^2 |u|_{2,K}^2. \tag{2.5}$$

**2. Estimate on irregular elements**.

Let $\widetilde{K}$ be any irregular element, which contains a subset $S_K \subset S$ (see Fig 1). Let $p_i = (x_i, y_i)$ be the vertices of $\widetilde{K}$. Since $S$ is piecewise smooth, there is a triangle $p_0 p_3 p_1$ in $\widetilde{K}$, with height $O(h^2)$, such that $S_K$ is contained in this triangle. Extend $p_1 p_0$ to intersect with $p_2 p_3$ at $p_1'$, and extend $p_3 p_0$ to intersect with $p_1 p_2$ at $p_3'$, respectively. Let $e$ denote the quadrilateral $p_1' p_0 p_3' p_2$, and $G$ denote the quadrilateral $p_0 p_1 p_2 p_3$. Since $u \in H^2(G)$, we employ an argument similar to the case of regular elements to obtain

$$\|u - u_I\|_{1,e}^2 \leqslant Ch^2 |u|_{2,G}^2. \tag{2.6}$$

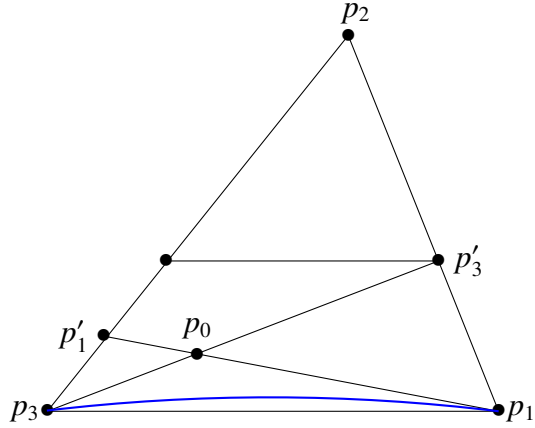

FIGURE 1. Illustration of an irregular cell $\widetilde{K}$.

Draw a trapezoid $T$ that has height $d$, which is the larger of the distances from $p'_1$ and $p'_3$ to edge $p_1 p_3$. Since $S_k$ passes through two vertices, it holds that $d = O(h^2)$. Now, we estimate $\|u - u_I\|_{1,T}$. In general, we only have $u \in H^1(T)$ on $T$. For this reason, we may only write

$$u(p) - u_I(p) = \sum_i L_i \int_0^1 \partial_t u(M_i) dt. \tag{2.7}$$

Again, we use the variable transformation $\zeta = x_i + (x - x_i)t$, $\eta = y_i + (y - y_i)t$, with $\dfrac{\partial(x,y)}{\partial(\zeta,\eta)} = t^{-2}$. Then, the region $T$ is changed to $T_{i,t}$, which is similar to $T$ with area $|T_{i,t}| \leq Ct^2 h^3$. For fixed $i, t$, let $D(i,t) = \bigcup T_{i,t}$ where the union is taken over all irregular elements. Since the number of irregular elements does not exceed $O(h^{-1})$, we have $|D(i,t)| \leq Ct^2 h^2$.

Given any $\varepsilon$ satisfying $0 < \varepsilon \leq \tfrac{1}{4}$, we have

$$\left| \int_0^1 \partial_t u(M_i) dt \right|^2 \leq \int_0^1 t^{-2\varepsilon} dt \int_0^1 t^{2\varepsilon} |\partial_t u(M_i)|^2 dt$$

$$\leq 2 \int_0^1 t^{2\varepsilon} |\partial_t u(M_i)|^2 dt.$$

Therefore,

$$\sum_T \iint_T \left| \int_0^1 \partial_t u(M_i) dt \right|^2 dxdy$$

$$\leq Ch^2 \int_0^1 t^{2\varepsilon - 2} \sum_T \iint_{T_{i,t}} \left[ u_x^2(\zeta, \eta) + u_y^2(\zeta, \eta) \right] d\zeta d\eta$$

$$\leq Ch^2 \int_0^1 t^{2\varepsilon - 2} \iint_{D(i,t)} (u_x^2 + u_y^2) dp.$$

Applying Lemma 2 on $D(i,t) \cap \Omega_1$ and $D(i,t) \cap \Omega_2$ separately, we obtain

$$\iint_{D(i,t)} (u_x^2 + u_y^2) dp \leq \frac{C}{\varepsilon} (th)^{2-2\varepsilon} \|u\|^2_{2,\Omega_1 + \Omega_2}.$$

Thus,
$$\sum_T \|u - u_I\|_{0,T}^2 \leq \frac{C}{\varepsilon} h^{4-2\varepsilon} \|u\|_{2,\Omega_1+\Omega_2}^2. \tag{2.8}$$

Now, we turn to estimating $\|\partial_x(u - u_I)\|_{0,T}$. Since $\sum_i L_i = 1$, it holds that $\sum_i \partial_x L_i = 0$. Hence, we have

$$\begin{aligned}\partial_x(u - u_I) &= \partial_x u - \sum_i \partial_x L_i u(p_i) \\ &= u_x(x,y) + \sum_i \partial_x L_i \left[u(x,y) - u(x_i, y_i)\right] \\ &= u_x(x,y) + \sum_i \partial_x L_i \int_0^1 \partial_t u(M_i) dt.\end{aligned}$$

Taking norm on both sides, we obtain

$$\|\partial_x(u - u_I)\|_{0,T} \leq \|u_x\|_{0,T} + \frac{C}{h} \left\|\int_0^1 \partial_t u(M_i) dt\right\|_{0,T}.$$

As shown above, we see
$$\sum_T \|u_x\|_{0,T}^2 \leq \frac{C}{\varepsilon} h^{2-2\varepsilon} \|u\|_{2,\Omega_1+\Omega_2}^2.$$

Combining the two inequalities above yields
$$\sum_T \|\partial_x(u - u_I)\|_{0,T}^2 \leq \frac{C}{\varepsilon} h^{2-2\varepsilon} \|u\|_{2,\Omega_1+\Omega_2}^2.$$

The estimate for $\partial_y(u - u_i)$ is similar.

Combining the two cases discussed above, we obtain
$$\|u - u_I\|_{1,\Omega_h} \leq \frac{C}{\sqrt{\varepsilon}} h^{1-\varepsilon} \|u\|_{2,\Omega_1+\Omega_2}. \tag{2.9}$$

Here, $C$ does not depend on $\varepsilon$, and $0 < \varepsilon < \frac{1}{4}$ is arbitrary. It is easy to show that for any fixed sufficiently small $h$, $\frac{1}{\sqrt{\varepsilon}} h^{1-\varepsilon}$ achieves its minimum when $\varepsilon = \frac{1}{2} |\ln h|^{-1}$. Plugging this $\varepsilon$ into (2.9), we obtain

$$\|u - u_I\|_{1,\Omega_h} \leq C h |\ln h|^{1/2} \|u\|_{2,\Omega_1+\Omega_2}$$

Finally, the inequality (1.7) follows from the basic inequality (1.6).

学生优秀论文

# 具有间断系数的二阶椭圆型方程的有限元解的敛速估计

数学系　许进超

## §1　引　言

对于系数适当光滑的二阶椭圆型方程的有限元解的收敛性问题已有较完善的结果[1]。但在工程实际中，常会遇到系数有间断的方程，用有限元解这类问题，其收敛性问题必须重新研究。目前这类问题的一维情形已有工作[5]，但对于更有理论和实际意义的多维问题的研究尚少见。冯康教授在文[3]中提到过间断线为折线的这一简单情形，但对于一般情形以及收敛性问题未见论及。笔者借毕业实习机会，拟就具有间断系数的二维 Dirichlet 问题的有限元（三角形线性元）解的收敛性作初步探讨。研究结果表明系数的间断对有限元解的收敛性并无多大影响

为书写简单起见，本文考虑如下模型问题

$$\begin{cases} -\partial_x(B\partial_x u)-\partial_y(B\partial_y u)+\sigma u=f(x,y) \\ \left(B\dfrac{\partial u}{\partial n}\right)_{1,s}=\left(B\dfrac{\partial u}{\partial n}\right)_{2,s} \\ u|_\Gamma=0 \end{cases} \tag{1.1}$$

这里 $\Omega$ 为平面有界连通区域，其边界 $\Gamma$ 分段光滑，分段向外凸，曲线 $S$ 分段光滑，将区域 $\Omega$ 分为 $\Omega_1, \Omega_2$ 两个子域，$B=B(x,y)$ 看作 $\Omega_i$ 上的函数时属于 $C^1(\overline{\Omega_i})$，$i=1,2$，$\sigma\in L^\infty(\Omega), f\in L^2(\Omega)$，且在 $\Omega$ 上满足

$$B(x,y)\geq B_0>0, \qquad \sigma\geq 0 \tag{1.2}$$

引入双线性型

$$a(u,v)=\iint_\Omega [B(\partial_x u\,\partial_x v+\partial_y u\,\partial_y v)+\sigma uv]dxdy \tag{1.3}$$

易证得下面的变分原理，

求解方程 (1.1) 等价于下述变分问题：求 $u\in H_0^1(\Omega)$，使对任意 $v\in H_0^1(\Omega)$ 成立

$$a(u,v)=(f,v) \tag{1.4}$$

由 (1.3) 定义的 $a(u,v)$ 仍是 $H_0^1(\Omega)$ 上连续、正定的双线性型。设变分问题 (1.4) 有解 $u\in H_0^1(\Omega,s)=\{u\mid u\in H_0^1(\Omega), u\in H^2(\Omega_1), u\in H^2(\Omega_2)\}$。它在有限维子空间 $S_h\subset H_0^1(\Omega)$ 中的近似解 $u_h$ 仍满足误差基本估计式：

$$\|u-u_h\|_{1,\Omega}\leq C\inf_{v\in S_h}\|u-v\|_{1,\Omega} \tag{1.5}$$

其中 $C$ 为常数（我们一律不加区别地用 $C$ 表示与 $h, u, v$ 无关的常数），特别地，若 $u_I \in S_h$ 是 $u$ 的插值函数，则有

$$\|u - u_h\|_{1,\Omega} \leq C \|u - u_I\|_{1,\Omega} \tag{1.6}$$

考虑 $\Omega$ 的三角剖分 $\Omega_h \subset \Omega$，假设其边界 $\Gamma_h$ 上的节点都在 $\Gamma$ 上，凡包含有间断线 $S$ 的单元必有两节点在 $S$ 上，$S$ 的不光滑点应作为节点，每单元包含有半径为 $ch$ 的圆（拟一致剖分）。$S_h$ 是在 $\Gamma_h$ 上为 $0$，在 $\Omega$ 上连续的分片线性的多项式函数类。显然 $S_h \subset H_0^1(\Omega)$。不包含间断线的单元称为正常单元，否则称为非正常单元。本文主要结果是：变分问题 (1.4) 在 $S_h$ 中的有限元解 $u_h$ 与真解 $u \in H^2(\Omega, S)$ 的偏差有估计

$$\|u - u_h\|_{1,\Omega_h} \leq Ch |\ln h|^{\frac{1}{2}} \|u\|_{2, \Omega_1 + \Omega_2} \tag{1.7}$$

## § 2    结论的证明

先给出两个引理，它们是 Sobolev 嵌入定理[6] 的精确化。

**引理 1**   设 $\Omega$ 为前述平面区域，$f \in L^2(\Omega)$，$p = (x, y)$，$Q = (\zeta, \eta)$，$B(p, Q)$ 是 $p, Q$ 的有界函数，当 $p \neq Q$ 时连续。定义

$$v(p) = \iint_\Omega \frac{B(p, Q)}{|p - Q|} f(Q) dQ$$

则有常数 $C$，使对任何可测集 $D \in \Omega$ 有

$$\iint_D v^2(p) dp \leq \frac{C}{\varepsilon} |D|^{1-\varepsilon} \|f\|_{L^2(\Omega)}^2$$

其中 $|D|$ 为 $D$ 的测度，$\varepsilon$ 为小于 $1$ 的任意正数。

证：记 $q = \frac{2}{\varepsilon} > 2$，由 Hölder 不等式得[6]

$$\iint_D |v|^q dp \leq C^q q^{\frac{q}{2}} \|f\|_{L^2(\Omega)}^{q-2} \iint_D dp \iint_\Omega f^2(Q) r^{-\frac{1}{2}} dQ$$

交换积分次序，即得

$$\iint_D |v|^q dp \leq C^q q^{\frac{q}{2}} \|f\|_{L^2(\Omega)}^q$$

再利用 Hölder 不等式，

$$\iint_D v^2 dp \leq |D|^{1-\frac{2}{q}} \|v\|_{L^q(D)}^2 \leq \frac{C}{\varepsilon} |D|^{1-\varepsilon} \|f\|_{L^2(\Omega)}^2$$

证毕。

利用此引理及 Sobolev 积分恒等式可得

**引理 2**   对 $H^1(\Omega)$ 中任意函数 $v$ 有

$$\iint_D v^2 dp \leq \frac{C}{\varepsilon} |D|^{1-\varepsilon} \|v\|_{1,\Omega}^2$$

其中 $D, \Omega, \varepsilon, C$ 如引理 1 所述

下面转向估计式 (1.7) 的证明

1.  正常单元上的估计

正常单元上的估计，自然可按熟知的Bramble-Hilbert引理① 得到，但我们采用的方法是利用变点的Taylor展开④，它是我们以后估计的基础。任取正常单元$K$，设 $p_i(x_i,y_i)$ 为其顶点，$p=(x,y)$ 为变点，$M_i=[x_i+(x-x_i)t, y_i+(y-y_i)t]$，$i=1,2,3$，$0 \leq t \leq 1$，在研究过程中不妨设 $v \in C^2(\Omega)$. 由带积分型余项的Taylor公式

$$u(p) - u(p_i) = (x-x_i)\partial_x u(p) + (y-y_i)\partial_y u(p) - \int_0^1 t\partial_i^2 u(M_i)dt$$

利用面积坐标$L_i$，可得线性插值函数 $u_I$ 的误差

$$u(p) - u_I(p) = -\sum_{i=1}^3 L_i \int_0^1 t\,\partial_i^2 u(M_i)dt \tag{2.1}$$

利用 Hölder 不等式得

$$\|u - u_I\|_{0,K}^2 \leq 3\sum_i \iint_K |\int_0^1 t\partial_i^2 u(M_i)dt|^2 dxdy$$

$$\leq 3\sum_i \iint_K \int_0^1 t^2 |\partial_i^2 u(M_i)|^2 dt\,dxdy$$

交换积分次序，对任意固定的 $t,i$，引进变量代换 $\zeta = x_i + (x-x_i)t, \eta = y_i + (y-y_i)t$，区域$K$变为子域 $K_{i,t} \subset K$. 于是得

$$\|u - u_I\|_{0,K}^2 \leq Ch^4 |u|_{2,K}^2 \tag{2.2}$$

又由(2.1)有

$$\partial_x(u - u_I) = -\sum_i \partial_x L_i \int_0^1 t\partial_i^2 u(M_i)dt - \sum_i L_i \int_0^1 t\partial_x \partial_i^2 u(M_i)dt$$

利用分部积分，后一积分为

$$\sum_i L_i \int_0^1 t\partial_x \partial_i^2 u(M_i)dt = \sum_i L_i \int_0^1 t[\partial_i^2 u_x(M_i)t + 2\partial_i u_x(M_i)]dt$$

$$= \sum_i L_i t^2 \partial_i u_x(M_i)|_0^1 = 0$$

注意到 $|\partial_x L_i| \leq C/h$，利用前面方法得

$$\|\partial_x(u - u_I)\|_{0,K}^2 \leq Ch^2 |u|_{2,K}^2 \tag{2.3}$$

同理 $$\|\partial_y(u - u_I)\|_{0,K}^2 \leq Ch^2 |u|_{2,K}^2 \tag{2.4}$$

综合以上三式得

$$\|u - u_I\|_{1,K}^2 \leq Ch^2 |u|_{2,K}^2 \tag{2.5}$$

2. 非正常单元上的估计

任取一非正常单元 $\widetilde{K}$ 如图，它包含 $S$ 的一段 $S_K$，$p_i(x_i,y_i)$ 为 $\widetilde{K}$ 的顶点，由 $S$ 的光滑性可在 $\widetilde{K}$ 内作一高为 $o(h^2)$ 的三角形 $p_0 p_2 p_1$ 使 $S_K$ 在此三角形内，延长两腰，分别交 $p_2 p_3$ 于 $p_1'$，交 $p_1 p_2$ 于 $p_3'$，记 $e$ 为四边形 $p_1' p_0 p_3' p_2$，$G$ 为四边形 $p_0 p_1 p_2 p_3$，由于 $u \in H^2(G)$，用类似方法可得

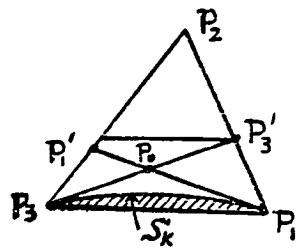

$$\|u - u_I\|_{1,e}^2 \leq Ch^2 |u|_{2,G}^2 \tag{2.6}$$

其次，以 $p_1', p_3'$ 到 $p_1 p_3$ 的较大距离 $d$ 为高作梯形 $T$，显然 $d = o(h^2)$. 我们估

计 $\|u-u_I\|_{1,T}$. 由于在 $T$ 上, 一般只有 $u \in H^1(T)$. 为此, 类似于 (2.1) 写出

$$u(p) - u_I(p) = \sum_i L_i \int_0^1 \partial_t u(M_i) dt \qquad (2.7)$$

仍作变量替换 $\zeta = x_i + (x-x_i)t, \eta = y_i + (y-y_i)t$, $\therefore \frac{\partial(x,y)}{\partial(\zeta,\eta)} = t^{-2}$, 区域 $T$ 变为 $T_{i,t}$, 且与 $T$ 相似, 面积 $|T_{i,t}| \leqslant C t^2 h^2$. 对固定的 $i,t$, 记 $D(i,t) = \cup T_{i,t}$, 求和是对所有非正常单元的, 由于非正常单元的个数不超过 $O(h^{-1})$, 故 $|D(i,t)| \leqslant C t^2 h$.

任取正数 $\varepsilon = > \frac{1}{4}$, 则

$$\left|\int_0^1 \partial_t u(M_i) dt\right|^2 \leqslant \int_0^1 t^{-2\varepsilon} dt \cdot \int_0^1 t^{2\varepsilon} |\partial_t u(M_i)|^2 dt$$
$$\leqslant 2 \int_0^1 t^{2\varepsilon} |\partial_t u(M_i)|^2 dt$$

故
$$\sum_T \iint_T |\int_0^1 \partial_t u(M_i) dt|^2 dx dy$$
$$\leqslant C h^2 \int_0^1 t^{2\varepsilon-2} \sum_T \iint_{T_{i,t}} [u_x^2(\zeta,\eta) + u_y^2(\zeta,\eta)] d\zeta d\eta$$
$$\leqslant C h^2 \int_0^1 t^{2\varepsilon-2} \iint_{D(i,t)} (u_x^2 + u_y^2) dp$$

分别在 $D(i,t) \cap \Omega_1$ 与 $D(i,t) \cap \Omega_2$ 上应用引理 2, 得

$$\iint_{D(i,t)} (u_x^2 + u_y^2) dp \leqslant \frac{C}{\varepsilon} (th)^{2-2\varepsilon} \|u\|_{2,\Omega_1+\Omega_2}^2$$

因此
$$\sum_T \|u-u_I\|_{0,T}^2 \leqslant \frac{C}{\varepsilon} h^{4-2\varepsilon} \|u\|_{2,\Omega_1+\Omega_2}^2 \qquad (2.8)$$

下面估计 $\|\partial_x(u-u_I)\|_{0,T}$. 因 $\sum_i L_i = 1$, 故 $\sum_i \partial_x L_i = 0$,

$$\partial_x(u-u_I) = \partial_x u - \sum_i \partial_x L_i u(p_i)$$
$$= u_x(x,y) - \sum_i \partial_x L_i [u(x,y) - u(x_i,y_i)]$$
$$= u_x(x,y) - \sum_i \partial_x L_i \int_0^1 \partial_t u(M_i) dt$$

取范数得

$$\|\partial_x(u-u_I)\|_{0,T} \leqslant \|u_x\|_{0,T} + \frac{C}{h} \|\int_0^1 \partial_t u(M_i) dt\|_{0,T}$$

与前面类似地有

$$\sum_T \|u_x\|_{0,T}^2 \leqslant \frac{C}{\varepsilon} h^{2-2\varepsilon} \|u\|_{2,\Omega_1+\Omega_2}^2$$

由以上两式得

$$\sum_T \|\partial_x(u-u_I)\|_{0,T}^2 \leqslant \frac{C}{\varepsilon} h^{2-2\varepsilon} \|u\|_{2,\Omega_1+\Omega_2}^2$$

估计 $\partial_y(u-u_I)$ 是类似的.

综合上面的两类估计, 我们得

$$\|u-u_I\|_{1,\Omega_h} \leqslant \frac{C}{\sqrt{\varepsilon}} h^{1-\varepsilon} \|u\|_{2,\Omega_1+\Omega_2} \qquad (2.9)$$

这里 $C$ 与 $\varepsilon$ 无关，$0<\varepsilon<\frac{1}{4}$ 是任意的. 容易证明，对固定的充分小的 $h$ ，取

$$\varepsilon=\frac{1}{2}|lnh|^{-1} \text{时} \frac{1}{\sqrt{\varepsilon}}h^{1-\varepsilon} \text{达到最小，将这个 } \varepsilon \text{ 代入(2.9)则有}$$

$$\|u-u_I\|_{1,\Omega_h} \leqslant Ch|lnh|^{\frac{1}{2}}\|u\|_{2,\Omega_1+\Omega_2}$$

最后由基本估计式 (1.6) 便得到估计式 (1.7).

### 参 考 文 献

# Estimate of the Convergence Rate of the Finite Element Solutions to Elliptic Equation of Second Order with Discontinuous Coefficients

Xu Jin—chao

## Abstract

In this paper, we consider elliptic boundary value problem with discontinuous coefficients and obtain the asymptotic optimal error estimate $\|u-u_h\|_{1,\Omega_h} \leqslant Ch|lnh|^{1/2}\|u\|_{2,\Omega_1+\Omega_2}$ for triangle linear elements.